\documentclass{article}
\usepackage{amsfonts}
\usepackage[mathscr]{eucal}
 \usepackage{amsmath}
 \usepackage{amssymb}
\usepackage{amsmath} 
  \usepackage{paralist}
  \usepackage{graphics} 
  \usepackage{epsfig} 
\usepackage{graphicx}  \usepackage{epstopdf}
 \usepackage[colorlinks=true]{hyperref}
\hypersetup{urlcolor=blue, citecolor=red}

\newtheorem{theorem}{Theorem}[section]

\newtheorem{proposition}{Proposition}

\newtheorem{definition}[theorem]{Definition}
\newtheorem{remark}{Remark}

\begin{document}

\begin{center}
{\Large  
Viral infection model with diffusion and state-dependent delay: stability of classical solutions\footnote{This paper is dedicated to the memory of Igor D. Chueshov}

\bigskip

Alexander  Rezounenko
}

\end{center} 

\centerline{V.N.Karazin Kharkiv National University,  Kharkiv, 61022, Ukraine}
   \centerline{Institute of Information Theory and Automation}
   \centerline{ Academy of Sciences  of the Czech Republic,  }
 \centerline{P.O. Box 18, 182\,08 Praha, CR }
 \centerline{ E-mail: rezounenko@gmail.com}









\bigskip


\begin{abstract}
A class of reaction-diffusion virus dynamics models with
intracellular state-dependent delay and
a general non-linear infection rate  functional response
is investigated. We are interested in classical solutions with Lipschitz in-time initial functions
which are adequate to the discontinuous change of parameters due to, for example, drug administration. The Lyapunov functions technique is used to analyse stability of interior infection equilibria which describe the cases of a chronic disease.
\end{abstract}



\def\R{\mathbb{R}}
\newcommand\dom{\operatorname{dom}}
\smallskip





\section{Introduction}

In our research we are interested in mathematical models  of viral
 diseases.
According to World Health Organization, many viruses (as  Ebola
virus, Zika virus, HIV, HBV, HCV and others) continue to be a major
global public health issues. Particularly, in the recent The Global hepatitis report (WHO, April 2017) we find  \cite{WHO-Global hepatitis report-2017}
 ''a large number of people - about 325 million
worldwide in 2015 - are carriers of hepatitis B or C virus infections, which can remain asymptomatic for decades.'' and ''Viral hepatitis caused 1.34 million deaths in 2015, a number comparable
to deaths caused by tuberculosis and higher than those caused by HIV.''
In such a situation any steps toward understanding viral diseases are important.

 There are variety of models with and without delays which describe
 dynamics of different viral infections. Delays could be
 concentrated or distributed, constant, time-dependent or state-dependent. 

We notice that classical models
\cite{Nowak-Bangham-S-1996,Perelson-Neumann-Markowitz-Leonard-Ho-S-1996}
contain ordinary differential equations (without delay) for  three variables: susceptible
host cells $T$, infected host cells $T^{*}$ and free virus particles $V$. The intracellular delay is an important property of the biological problem, so we formulate the delay problem
\begin{equation}\label{sdd-vir3-02}
\left\{
\begin{array}{l}
   \dot T(t) = \lambda - d T(t) - f(T(t),V(t)), \\
 \dot T^{*}(t) = e^{-\omega h} f(T(t-h),V (t-h))
 - \delta T^{*}(t), \\
  \dot V(t) = N\delta T^{*}(t) - cV (t). \\
\end{array}
\right.
\end{equation}
In (\ref{sdd-vir3-02}), susceptible cells $T$  are produced at a rate $\lambda$, die at rate $d T$, and become infected at rate $f(T,V)$. Properties and examples of incidence function $f$ are discussed below. Infected cells $T^{*}$  die at rate $\delta T^{*}$, free virions $V$ are produced by infected cells at rate  $ N\delta T^{*}$ and are removed at rate $cV (t)$.
In (\ref{sdd-vir3-02}) $h$ denotes the delay between the time a
virus particle contacts a target cell and the time the cell becomes
actively infected (start producing new virions). It is clear that
the constancy of the delay is an extra assumption which essentially
simplifies the analysis, but has no biological background.

 To the best of our knowledge, viral infection models with state-dependent delay (SDD)
 have  been considered for the first time in \cite{Rezounenko-DCDS-B-2017} (see also \cite{Rezounenko-EJQTDE-2016}).
It is well known that differential equations with state dependent
delay are always non-linear by its nature (see the review
\cite{Hartung-Krisztin-Walther-Wu-2006} for more details and
discussion).

  As usual in a delay system with
 (maximal) delay $h>0$ \cite{Hale,Kuang-1993_book,Walther_book}, for a function $v(t), t\in [a-h,b]\subset \mathbb{R}, b>a$,
 we denote the history segment $v_t=v_t(\theta)\equiv v(t+\theta), \theta\in [-h,0], t\in [a,b].$

The ODEs delay system (\ref{sdd-vir3-02}) is extended to the state-dependent one
\begin{equation}\label{sdd-vir3-03}
\left\{
\begin{array}{l}
   \dot T(t) = \lambda - d T(t) - f(T(t),V(t)), \\
 \dot T^{*}(t) = e^{-\omega h} f(T(t-\eta(u_t)),V (t-\eta(u_t))) - \delta T^{*}(t), \\
  \dot V(t) = N\delta T^{*}(t) - cV (t). \\
\end{array}
\right.
\end{equation}
Here $u(t)=(T(t),T^{*}(t),V(t))$.
System (\ref{sdd-vir3-03}) is a particular case of the system with state-dependent delay studied in \cite{Rezounenko-DCDS-B-2017,Rezounenko-EJQTDE-2016}. 
The ODE system is formulated assuming  host cells do not move and
the diffusion of free virus particles is very quick, so  they are
mixed enough to consider homogeneous distribution over the spatial
domain in a host organ. Similar situation is in case of all cells
and free virions are well mixed (e.g., in case of HIV and other
infections targeting blood cells). To consider more realistic
nonhomogeneous  situation one introduces spatial coordinate $x\in \Omega$ and allow the unknowns to depend on it, i.e.
$T(t,x),T^{*}(t,x),V(t,x)$.
Now $T(t,x),T^{*}(t,x),V(t,x)$ represent the densities of uninfected cells, infected cells and free virions at position $x$ at time $t$.

Consider a
connected
bounded domain $\Omega \subset \mathbb{R}^n$ with a smooth boundary $\partial \Omega$.
 Now we are
 ready to  present the PDEs system under consideration
\begin{equation}\label{sdd-vir3-04}
\left\{
\begin{array}{l}
   \dot T(t,x) = \lambda - d T(t,x) - f(T(t,x),V(t,x)) + d^1 \Delta T(t,x), \\
 \dot T^{*}(t,x) = e^{-\omega h} f(T(t-\eta(u_t),x),V (t-\eta(u_t),x)) - \delta T^{*}(t,x) + d^2 \Delta T^{*}(t,x), \\
  \dot V(t,x) = N\delta T^{*}(t,x) - cV (t,x) + d^3 \Delta V(t,x). \\
\end{array}
\right.
\end{equation}
Here the dot over a function denotes the partial time derivative i.g, $\dot
T(t,x) = {\partial T(t,x)\over \partial t}$, all the constants $\lambda, d, \delta,
 N, c,\omega$ are positive while $d^i, i=1,2,3$ (diffusion coefficients) are non negative.
We consider a general functional response $f(T,V)$ satisfying natural assumptions presented below. In earlier models (with constant or without delay) the study was started in case of bilinear $f(T,V) =\mbox{const}\cdot T V$ and then extended to more general classes of non-linearities, see Remark \ref{sdd-vir3-rem1} below.

Boundary conditions are of Neumann type for the corresponding unknown if $d^i\neq 0$ i.e.
${\partial T(t,x)\over \partial n} |_{\partial \Omega}=0$ if $d^1\neq 0$
and similarly for $T^{*}(t,x)$ and $V(t,x)$. Here ${\partial \over \partial n}$
is the outward normal derivative on $\partial \Omega$.
In case $d^i= 0$, no boundary conditions are needed for the corresponding unknown(s).

 Our main goals are to present the existence and uniqueness results for the  model (\ref{sdd-vir3-04}) in the sense of  classical solutions, and to study the local asymptotic stability of  non-trivial diseased equilibria. We apply the Lyapunov approach \cite{Lyapunov-1892} to the state-dependent delay PDE model and allow, but not require, diffusion terms in each state equation.

There is a number of works studying the case $d^1=d^2=0, d^3>0$  (see e.g. \cite{Wang-Wang-MB-2007_HBV_spatial dependence,Wang-Huang-Zou-AA-2014,Wang-Yang-Kuniya-JMAA-2016} for models without delay and
\cite{McCluskey-Yang-NA-2015,Hattaf-Yousfi-CMA-2015} with {\it constant} delay; see also references therein). In the mentioned works  authors assume that the host cells (healthy and infected) do not move or are well mixed, while viral particles diffuse freely. Let us discuss the cases when an infection affects one particular organ as, for example, liver in case of HBV, HCV. In
such cases the spatial domain $\Omega \subset \mathbb{R}^3$
represents the organ. The Neumann boundary conditions say that viral particles do not leave the organ. It is not relevant from the biological point of view since viral particles circulate together with the blood stream in and out the organ (e.g. liver). For the mathematical system to cover such cases one could assume $d^3=0$ and
no boundary conditions for $V$. Taking into account the high speed of the blood stream, this means the viral particles are well mixed.
Even more interesting case is $d^i>0, i=1,2, d^3=0$. To the best of our knowledge, this case has not been considered before. The case
$d^2>0$ may reflect the cell-to-cell transmission  of the
infection when viral particles cross the membranes of the nearest cells (see \cite{Carloni et al-CMM-2012} for more discussion and references; c.f. \cite{Wang-Yang-Kuniya-JMAA-2016}). The infection spreads similar to diffusion to cells in a neighbourhood of an infected cell. The case $d^1>0$ may reflect natural division of healthy cells in order to fill the space previously occupied by infected cells (after the death of the last ones). In cases $d^1>0, d^2>0$, the host cells (both healthy and infected) do not leave the organ, so Neumann boundary conditions are
quite relevant.

In study of state-dependent delay equations the choice of the set of initial functions is particularly important and non-trivial (see review \cite{Hartung-Krisztin-Walther-Wu-2006} for ODE case  and works \cite{Rezounenko_JMAA-2007,Rezounenko_NA-2009,Rezounenko_NA-2010,Chueshov-Rezounenko_CPAA-2015} for PDEs).
%
%
We are interested in classical solutions with Lipschitz in-time initial functions which are adequate to the discontinuous change of parameters due to, for example, drug administration (for more discussion and references see \cite{Rezounenko-EJQTDE-2016}).
The main motivation here is the situation (see e.g.
\cite{Shudo-Ribeiro-Talal-Perelson_Antiviral
Therapy-2008,Pawlotsky-Semin Liver Dis-2014}) when the drug
effectiveness is decreased in a stepwise manner. In terms of
system (\ref{sdd-vir3-04}), the parameter $N$ could change its value in a discontinuous way (see equation (2) in
\cite[p.920]{Shudo-Ribeiro-Talal-Perelson_Antiviral Therapy-2008}).
It is clear that at any time moment of discontinuity of (any) parameter, the solution is continuous, but not differentiable (c.f. figure 2-B in \cite[p.921]{Shudo-Ribeiro-Talal-Perelson_Antiviral
Therapy-2008} and also fig.1 in \cite[p.23]{Pawlotsky-Semin Liver Dis-2014}).

Since {\it delay} is a central part of the paper, it would be
interesting to present examples of SDD $\eta$ and discuss the structure of $\eta$ from biological point of view. Unfortunately, up to now, the biological side of virus dynamics  is not fully understood. Even current {\it in vitro} study does not provide enough information. {\it In vivo} study is essentially more complicated, and up to now, there are no technical (biological/medical) tools for the {\it real time} monitoring of disease dynamics available.
In such a situation we present a rather general class of SDD (see (\ref{sdd-vir3-example1}),  (\ref{sdd-vir3-example2}) below). Delays of the form (\ref{sdd-vir3-example1}), (\ref{sdd-vir3-example2}) take into account all the prehistory $u_t$ by integrating a solution over $[t-h,t]$.

For general facts on PDEs with {\it constant delay} see e.g.
\cite{travis_webb,Martin-Smith-TAMS-1990,Wu_book} and PDEs with {\it
state-dependent delay} 
\cite{Rezounenko_JMAA-2007,Rezounenko_NA-2009,Rezounenko_NA-2010,Rezounenko_JMAA-2012,Rezounenko-JADEA-2012,Rezounenko-Zagalak-DCDS-2013,Chueshov-Rezounenko_CPAA-2015}.
We also mention that the case of all $d^i>0$ is, in a sense, easier
from mathematical point of view since the linear part generates a
compact semi-group.

We use the Lyapunov functions technique \cite{Lyapunov-1892} to analyse stability of interior infection equilibria which describe the cases of chronic disease.
To the best of our knowledge, viral infection models with diffusion
and state-dependent delay have not been considered before.

\section{Basic properties of the model}

Define the following linear operator \ $ -\mathcal{A}^0 $\ $= diag \, (d^1 \Delta, d^2 \Delta,
d^3 \Delta)$ in $C(\overline{\Omega}; \mathbb{R}^3)$ with
$D(\mathcal{A}^0)\equiv D(d^1 \Delta)\times D(d^2 \Delta) \times D(d^3
\Delta)$. Here, for $d^i\neq 0$ we set $D(d^i \Delta)\equiv \{ v\in
C^2(\overline{\Omega}) : {\partial v(x)\over \partial n} |_{\partial
\Omega}=0\}$
 and $D(d^j \Delta)\equiv C(\overline{\Omega})$ for $d^j= 0$.
We omit the space coordinate $x$, for short, for unknown
$u(t)=(T(t),T^{*}(t),V(t))\in X\equiv
[C(\overline{\Omega})]^3\equiv  C(\overline{\Omega}; \mathbb{R}^3)$.
It is well-known that
the closure $-\mathcal{A}$ (in $X$) of
the  operator $-\mathcal{A}^0$ generates a
$C_0$-semigroup $e^{-\mathcal{A} t}$ on $X$ which is analytic and
nonexpansive \cite[p.5]{Martin-Smith-TAMS-1990}. We denote the space
of continuous functions by $C\equiv C([-h,0]; X)$ equipped with the
sup-norm $||\psi||_C\equiv \max_{\theta\in [-h,0]} ||\psi
(\theta)||_X$.

 We  write,   the system (\ref{sdd-vir3-04}) in abstract form
 \begin{equation}\label{sdd-vir3-abst-eq}
 {d\over dt} u(t) + \mathcal{A} u(t) = F(u_t), \qquad t>0.
 \end{equation}
 The non-linear continuous mapping $F: C \to X$ is defined by
 \begin{equation}\label{sdd-vir3-F}
 F(\varphi)=  F(\varphi)(x)=\left(\begin{array}{l}
  \lambda - d \varphi^1(t,x) - f(\varphi^1(t,x),\varphi^3(t,x))  \\
 e^{-\omega h} f(\varphi^1(-\eta(\varphi),x),\varphi^3 (-\eta(\varphi),x)) - \delta \varphi^2(t,x)  \\
   N\delta \varphi^2(t,x) - c \varphi^3 (t,x)\\
 \end{array}
 \right).
 \end{equation}
Here $ \varphi = ( \varphi^1,\varphi^2,\varphi^3)\in C$.  Mapping $F$ is {\it not} Lipschitz on the space $C$ which is typical
 for a mapping which includes discrete state-dependent delays
 (see review \cite{Hartung-Krisztin-Walther-Wu-2006} for ODE case
 and works \cite{Rezounenko_JMAA-2007,Rezounenko_NA-2009,Rezounenko_NA-2010,Chueshov-Rezounenko_CPAA-2015} for PDEs).

We need initial conditions $u(\theta,x) =\varphi (\theta,x) = (T(\theta,x),T^{*}(\theta,x),V(\theta,x)),  \theta\in [-h,0]$ for the delay problem (\ref{sdd-vir3-abst-eq})

\begin{equation}\label{sdd-vir3-ic2}
\varphi \in Lip ([-h,0]; X)\equiv \left\{
\psi\in C : \sup_{s\neq t} {||\psi(s)-\psi(t)||_X\over |s-t|} < \infty
\right\}, \quad
\varphi(0)\in D(\mathcal{A}).
\end{equation}
In our study we use the standard (c.f. \cite[Def. 2.3,
p.106]{Pazy-1983-book} and \cite[Def. 2.1, p.105]{Pazy-1983-book})
\smallskip

\begin{definition}\label{def1}
 A function $u\in C([-h,T]; X)$ is called a {\tt mild solution}
 on $[-h,T)$ of the initial value problem (\ref{sdd-vir3-abst-eq}),
 (\ref{sdd-vir3-ic2}) if it satisfies
  (\ref{sdd-vir3-ic2}) and
 \begin{equation}\label{sdd-vir3-3-1}
u(t)=e^{-\mathcal{A} t}\varphi(0) + \int^{t}_0 e^{- \mathcal{A} (t-s)}  F(u_s) \, ds,
\quad t\in [0,T).
 \end{equation}
A function $u\in C([-h,T); X)\bigcap C^1((0,T);
X)$ is called a {\tt classical solution}  on $[-h,T)$ of the initial value problem (\ref{sdd-vir3-abst-eq}), (\ref{sdd-vir3-ic2}) if it
satisfies  (\ref{sdd-vir3-ic2}), $u(t)\in D(A)$ for $0<t<T$ and
(\ref{sdd-vir3-abst-eq}) is satisfied on $(0,T)$.
\end{definition}

In the study below we are mainly interested in classical solutions which preserve the regularity of the Lipschitzian initial data (see (\ref{sdd-vir3-ic2})).

\medskip


Assume the non-linear function $f: \mathbb{R}^2\to \mathbb{R}$ is Lipschitz continuous and satisfies
 \begin{equation}\label{sdd-vir3-f-1}
 {\bf (Hf_1)} \quad
  \text{there exists } \mu >0 \,
  \text{ such that  } |f(T,V)|\le \mu |T| \,\, \text{for all } T, V\in \mathbb{R},
 \end{equation}

We have the following result

\begin{proposition}\label{proposition1}
Let nonlinear function $f$ be Lipschitz and  satisfy
 ${\bf (Hf_1)}$ (see (\ref{sdd-vir3-f-1})),
state-dependent delay $\eta : C \to [0,h]$ is locally Lipschitz.
Then
the initial value problem (\ref{sdd-vir3-abst-eq}), (\ref{sdd-vir3-ic2})
has a unique classical 
solution which is global in time i.e. defined for all $t\ge 0$.
\end{proposition}

{\it Proof of Proposition \ref{proposition1}}. We start with
discussion of mild solutions. Since the semigroup  generated by the
linear part $-\mathcal{A}$ is not necessarily compact (in cases when at
least one constant $d^i\neq 0$), see e.g.
\cite{Martin-Smith-TAMS-1990}, we cannot directly use results of
\cite{Rezounenko_JMAA-2007,Rezounenko_NA-2009,Rezounenko_NA-2010,Chueshov-Rezounenko_CPAA-2015}.
On the other hand, as mentioned above, non-linearity $F$ is not
Lipschitz on $C$, so we cannot directly apply the existence result
of \cite{Martin-Smith-TAMS-1990}. Moreover, the extension provided
in \cite{Rezounenko-JADEA-2012} cannot be directly applied to our
case since we do not assume here the ignoring condition on the
state-dependent delay (see more details in \cite{Rezounenko_NA-2009,
Rezounenko_JMAA-2012,Rezounenko-JADEA-2012}). Nevertheless, the
restrictions on initial function $\varphi$ posed by
(\ref{sdd-vir3-ic2}) give the possibility  to prove the existence of
a (unique) mild solution to initial-value problem
(\ref{sdd-vir3-04}), (\ref{sdd-vir3-ic2}) using the standard line
based on Banach Fixed Point Theorem  (in a complete metric space) as in the ODE case.
We outline only main steps of the proof. First we consider the following extension of $\bar\varphi (t)=\varphi (t)$
for $t\in [-h,0]$ and $\bar\varphi (t)=e^{-\mathcal{A} t}\varphi (0)$
for $t\ge 0$. Next we change variable $u(t)=\bar\varphi (t) +y(t)$
and consider complete metric space $A(\alpha,\beta,\gamma)\equiv \{
y\in C([-h,\alpha];X), y_0\equiv 0, \max_{t\in [0,\alpha]}
||y(t)||_X\le\beta, \sup_{s\neq t} ||y(s)-y(t)||_X\cdot |s-t|^{-1}
\le \gamma\}$ endowed by the metrics of the space of continuous
functions.  The operator $\mathcal{F}: A(\alpha,\beta,\gamma) \to
C([-h,\alpha];X)$ is defined as $\mathcal{F}(y)(t)\equiv 0$ for
$t\in [-h,0]$ and $\mathcal{F}(y)(t)\equiv \int^t_0 e^{\mathcal{A}
(t-\tau)} F(\bar\varphi_\tau + y_\tau) \, d\tau$ for $t\in
(0,\alpha]$. It is not difficult to check that our non-linear mapping
F, defined by (\ref{sdd-vir3-F}), satisfies (see the estimate $|f(T,V)|\le
\mu |T|$ in ${\bf (Hf_1)}$) $||F(\psi)||_X\le n_1 + n_2 ||\psi||_C$
and is locally {\it almost Lipschitz} on $A(\alpha,\beta,\gamma)$ by
the terminology of \cite{Mallet-Paret}. The last means
$||F(\psi^1)-F(\psi^2)||_X \le L_F(\gamma) ||\psi^1-\psi^2||_C$.
Standard computations show that operator $\mathcal{F}$ maps
$A(\alpha,\beta,\gamma)$ into itself provided $\alpha, \beta,\gamma$
satisfy $\alpha (n_1+n_2 (||\varphi|| + \beta))\le \beta, n_1+n_2
\beta \le \gamma$. Additional condition $\alpha L_F(\gamma)  < 1$
guarantees the contraction of $\mathcal{F}$. The classical Banach
Fixed Point Theorem gives the unique fixed point $\hat{y}$ and hence
the unique mild solution $u=\bar\varphi + \hat{y}$. The linear
growth bound of $F$ implies the global continuation of the mild
solution.

Our next step is to show that any mild solution is classical.
Let us fix any mild solution $u$ to (\ref{sdd-vir3-abst-eq}), (\ref{sdd-vir3-ic2}) and define $g(t)\equiv F(u_t), t\ge 0.$
For any $t^0>0$, mapping $g$  is continuous on $[0,t^0]$ since $F$ and $u$ are continuous.
We notice that by construction, the solution is Lipschitz
in time on $[0,t^0]$ (see also restrictions in
(\ref{sdd-vir3-ic2})). Hence,
$||g(t)-g(s)|| = ||F(u_t)-F(u_s)|| \le L_F \max_{\theta \in [-h,0]} ||u(t+\theta)-u(s+\theta)|| \le L_F L^{[0,t^0]}_u \cdot |t-s|$. Here we use the {\it almost Lipschitz} property of $F$.
Now we consider the following (non-delayed) initial value problem
\begin{equation}\label{sdd-vir3-non-delayed}
{d v(t)\over dt} + \mathcal{A} v(t) = g(t), \quad v(0)=x \in X,
\end{equation}
which has a unique solution. The solution of (\ref{sdd-vir3-non-delayed}) is $v=u$ in case $x=u(0)$.

We remind that $C_0$-semigroup
$e^{-\mathcal{A} t}$ is analytic on $X$
 \cite[p.5]{Martin-Smith-TAMS-1990}.   Hence theorem 3.5
 \cite[p.114]{Pazy-1983-book} implies that the mild solution (of  (\ref{sdd-vir3-non-delayed}) and hence of (\ref{sdd-vir3-abst-eq}), (\ref{sdd-vir3-ic2})) is classical for $t\ge 0.$
 The proof of Proposition \ref{proposition1} is complete.


Define the set (c.f. (\ref{sdd-vir3-ic2}))
$$
\Omega_{Lip} \equiv \left\{
\varphi = (\varphi^1,\varphi^2,\varphi^3)
\in Lip ([-h,0]; X))\subset C, \, \, \varphi(0)\in D(\mathcal{A})
 : \quad 0\le \varphi^1(\theta)\le {\lambda\over d},
\right.
$$
\begin{equation}\label{sdd-vir3-omega}
\left.
0\le \varphi^2(\theta)\le {\lambda\mu\over d\delta} e^{-\omega h}, \quad 0\le \varphi^3(\theta)\le {N\lambda\mu\over d c} e^{-\omega h}, \quad \theta \in [-h,0]
\right\},
\end{equation}
where $\mu$ is defined in $(Hf_1)$
and all the inequalities hold pointwise w.r.t. $x\in\overline{\Omega}.$

 We need further assumptions (which include $(Hf_1)$) on Lipschitz 
 function $f$ :
 \begin{equation}\label{sdd-vir3-f-1+}
 {\bf (Hf_1+)} \quad \left\{
 \begin{array}{l}
f(T,0)=f(0,V)=0, \quad \text{ and }\quad  f(T,V)>0 \text{ for all } T>0,V>0;  \\
 f \text{ is strictly increasing in both coordinates for all } T>0,V>0;   \\
  \text{there exists } \mu >0 \,
  \text{ such that  } |f(T,V)|\le \mu |T| \,\, \text{for all } T, V\in \mathbb{R}.
 \end{array}
  \right.
 \end{equation}



 We have the following result

\begin{proposition}\label{proposition2}
Let non-linear function $f$ satisfy $(Hf_1+)$ (see
(\ref{sdd-vir3-f-1+})),
state-depen-dent delay $\eta : C \to [0,h]$ is locally Lipschitz.
Then $ \Omega_{Lip}$ is invariant i.e. for any $\varphi \in
\Omega_{Lip}$ the unique solution to problem
(\ref{sdd-vir3-abst-eq}), (\ref{sdd-vir3-ic2}) satisfies $u_t \in
\Omega_{Lip}$  for all $t\ge 0$.
\end{proposition}

{\it Proof of Proposition \ref{proposition2}}. The existence and
uniqueness of solution is proven in Proposition \ref{proposition1}.
The proof of the  invariance part follows  the invariance result
of \cite{Martin-Smith-TAMS-1990} with the use of the almost
Lipschitz property of nonlinearity $F$. The estimates (for  the
subtangential condition) are the same as for  the constant delay
case, see e.g. \cite[Theorem 2.2]{McCluskey-Yang-NA-2015}. We do not
repeat it here. It is important to notice that the solutions are
classic for all $t\ge 0$ (but not for $t\ge h$ as could be in the
case of merely continuous initial functions $\varphi\in C$). The
proof of Proposition \ref{proposition2} is complete. 

\subsection{Stationary solutions}

Let us discuss stationary solutions of (\ref{sdd-vir3-04}).
By such solutions we mean time independent $\widehat u$ which,
in general, may depend on $x\in\overline{\Omega}$.
 Consider the system (\ref{sdd-vir3-04}) with
$u(t)=u(t-\eta(u_t))=\widehat u$ and denote the coordinates of a
stationary solution by
$(\widehat{T},\widehat{T^{*}},\widehat{V})=\widehat
u\equiv \widehat{\varphi}(\theta),\, \theta\in [-h,0]$.
Since stationary solutions of (\ref{sdd-vir3-04}) do not depend
on the type of delay (state-dependent or constant)
we have (see e.g. \cite{McCluskey-Yang-NA-2015})
\begin{equation}\label{sdd-vir3-stationary-1}
 \left\{
\begin{array}{l}
   0 = \lambda - d \widehat{T} - f(\widehat{T},\widehat{V}), \qquad
 0 = e^{-\omega h} f(\widehat{T},\widehat{V}) - \delta \widehat{T^{*}}, \\
  0 = N\delta \widehat{T^{*}} - c \widehat{V}. \\
\end{array}
\right.
\end{equation}
Equations  hold pointwise w.r.t. $x\in\overline{\Omega}.$

It is easy to see that the trivial stationary solution $(\lambda d^{-1},0,0)$ always exists.
We are interested in nontrivial disease stationary solutions of
(\ref{sdd-vir3-04}).
Using (\ref{sdd-vir3-stationary-1}), we have $\widehat{T}=(\lambda -
\delta \widehat{T^{*}}e^{\omega h})d^{-1} $ and
$\widehat{V}={N\delta\over c} \widehat{T^{*}}$. It gives the
condition on the coordinate $\widehat{T^{*}}$ which should belong to
$(0,\lambda e^{\omega h}\delta^{-1}]$. Denote (c.f.
\cite{McCluskey-Yang-NA-2015})
\begin{equation}\label{sdd-vir3-stationary-2}
h_f(s)\equiv f \left( {\lambda \over d}- {\delta \over d}e^{\omega
h}\cdot s,{N\delta\over c} \cdot s \right) - \delta e^{\omega
h}\cdot s.
\end{equation}
Assume $f$ satisfies
\par\medskip
${\bf (Hf_2)} \hskip10mm h_f(s)=0$ has at least one and at most finite roots on $(0,\lambda e^{\omega h}\delta^{-1}]$.
\medskip


We denote an arbitrary root of $h_f(s)=0$ by $\widehat{T^{*}}$ and
define the corresponding $\widehat{T}=(\lambda - \delta
\widehat{T^{*}}e^{\omega h})d^{-1} $ and $\widehat{V}={N\delta\over
c} \widehat{T^{*}}$. The point $(\widehat{T},
\widehat{T^{*}},\widehat V)$ satisfies
(\ref{sdd-vir3-stationary-1}), so it is a disease stationary
solutions of (\ref{sdd-vir3-04}).

\begin{remark}\label{sdd-vir3-rem-stationary}
We notice that the finiteness  of roots (which are obviously  isolated) does not allow the existence of equilibria which depend on spatial coordinate $x\in \Omega$. We remind that $\Omega$ is a connected set, so a function $v\in C(\overline{\Omega})$  may take either one or continuum values. Assumption  ${\bf (Hf_2)}$ implies $\widehat{T^{*}}(x)\equiv \widehat{T^{*}}\in \mathbb{R}$, so $(\widehat{T},
\widehat{T^{*}},\widehat V)$ is independent of  $x\in\overline{\Omega}.$
\end{remark}

\begin{remark}\label{sdd-vir3-rem1}
Below we mention some well-known examples of non-linear functions
$f$ when we have exactly one root of $h_f(s)=0$.
The first one is the DeAngelis-Bendington
\cite{Beddington-JAE-1975,DeAngelis-Goldstein-ONeill-E-1975} functional response
 $f(T,V)={kTV\over 1+k_1 T+k_2 V}$, with $k,k_1\ge
0,k_2>0$. We also mention
that the  functional response includes as a
special case ($k_1=0$) the {\it saturated incidence} rate
$f(T,V)={kTV\over 1+k_2 V}$.
 Another example of the nonlinearity is the Crowley-Martin
incidence rate $f(T,V)={kTV\over (1+k_1 T)(1+k_2 V)}$, with $k\ge
0,k_1,k_2>0$ (see e.g. \cite{Xu_JQTDE-2012}).
 For more general class of functions $f$ see, e.g.
 \cite{McCluskey-Yang-NA-2015,Hattaf-Yousfi-CMA-2015,Rezounenko-EJQTDE-2016}, where  under
additional conditions, one has exactly one root of $h_f(s)=0$. We notice that, in contrast
to \cite{McCluskey-Yang-NA-2015,Hattaf-Yousfi-CMA-2015}, we do not assume here
the differentiability of $f$.
\end{remark}
\begin{remark}\label{sdd-vir3-v-functional-3}
It is important to mention that usually in study of stability properties of stationary solutions (for viral dynamics problems) one uses conditions on the so-called  reproduction numbers. These conditions are used to separate the case of a unique stationary solution. Then the global stability of the equilibrium is investigated. In our study, taking into account the state-dependence of the delay, we discuss the  local stability. As a consequence, it allows  the co-existence of multiple equilibria. We believe this framework provides a way to model more complicated situations with rich dynamics (in contrast to a globally stable equilibrium). The conditions on the reproduction numbers do not appear explicitly here, but could be seen as particular sufficient conditions for ${\bf (Hf_2)}$.
\end{remark}


\section{Stability of disease stationary solutions}

The following Volterra function $v(s) = s - 1- \ln s: (0,+\infty)
\to \mathbb{R}_{+}$ plays an important role in construction of
Lyapunov functionals
\cite{Korobeinikov-BMB-2007,McCluskey-Yang-NA-2015}. One can see
that $v(s)\ge 0$ and $v(s)=0$ if and only if $s=1$. The derivative
equals $\dot v(s) = 1-{1\over s}$, which is obviously negative for
$x\in (0,1)$ and positive for $x>1$. The graph of $v$ explains the
use of the composition $v\left({s\over s^0}\right)$ in the study of
the stability properties of an equilibrium $s^0$. Another important
property is the following  \cite{Rezounenko-DCDS-B-2017} estimate
\begin{equation}\label{sdd-vir3-v-functional-2}
\forall \mu \in (0,1)\quad \forall s \in (1-\mu, 1+\mu) \quad \mbox{
one has } \quad {(s-1)^2\over 2(1+\mu)} \le v(s) \le {(s-1)^2\over
2(1-\mu)}.
\end{equation}
To check it, one simply observes that all three functions vanish at
$s=1$ and $\left|{d\over ds} \left( {(s-1)^2\over 2(1+\mu)} \right)
\right| \le  |{d\over ds}  v(s)| \le \left|{d\over ds} \left(
{(s-1)^2\over 2(1-\mu)} \right) \right|$ in the $\mu$-neighbourhood
of $s=1$.


In this section we use the following {\it local} assumptions
on $f$ in a small neighbourhood of a disease equilibrium
(given by ${\bf (Hf_2)}$).

 \begin{equation}\label{sdd-vir3-12}
{\bf (Hf_3)} \hskip30mm \left(   {V\over \widehat{V }} -
 {f(T,V)\over f(T,\widehat{V}) } \right) \cdot \left( {f(T,V)\over f(T,\widehat{V}) } -1\right) > 0.
\end{equation}
  This property simply means that the value ${f(T,V)\over f(T,\widehat{V}) }$ is always {\it strictly} between $1$ and ${V\over \widehat{V }}$ for any $T\ge 0$ (c.f. with the non-strict property \cite[p.74]{McCluskey-Yang-NA-2015}). The strict inequality in (\ref{sdd-vir3-12}) will be needed to handle the state-dependence of the delay. In the particular case of constant delay, the  non-strict property is enough.

We will also use the following assumption

\smallskip

${\bf (Hf_4)}$ Function $f$ is either differentiable with respect
to its first coordinate or
satisfies
 \begin{equation}\label{sdd-vir3-f-3}
 [f(T,\widehat{V})]^{-1} \ge C^1_f + C^2_f {1\over T},\quad T>0,
 \quad C^i_f = C^i_f(\widehat{V})\ge 0,\, i=1,2.
 \end{equation}

\smallskip
For simplicity of presentation we start with stability analysis
for smooth initial data belonging to the so-called {\it solution manifold}
(see e.g. \cite{Walther_JDE-2003,Hartung-Krisztin-Walther-Wu-2006}
for ODE case and \cite{Rezounenko-Zagalak-DCDS-2013} for PDEs)
 \begin{equation}\label{sdd-vir3-M_F}
M_F\equiv \left\{
\varphi \in C^1([-h,0];X),\quad \varphi(0)\in D(\mathcal{A}),
\quad \dot \varphi (0)+ \mathcal{A}\varphi (0) = F(\varphi) \,  \right\} .
 \end{equation}

The equation in (\ref{sdd-vir3-M_F}), called the compatibility
conditions, is an equality in $X$. Below (see Theorem \ref{vir3-theorem-stab2})
we return to more general case of  Lipschitz initial functions
($\varphi \in \Omega_{Lip}$,  not necessarily continuously differentiable)
which are important to cover the cases of drug administration when  the time derivative may be discontinuous, see
 \cite{Rezounenko-EJQTDE-2016} for more discussion.

\begin{theorem}\label{vir3-theorem-stab1}
Let the nonlinear function $f$ satisfy $(Hf_1+), (Hf_2), (Hf_3), (Hf_4)$
(see (\ref{sdd-vir3-f-1+}), (\ref{sdd-vir3-f-3}),
(\ref{sdd-vir3-12})) and  state-dependent delay $\eta : C \to [0,h]$
be locally Lipschitz in $C$ and continuously differentiable
in a neighbourhood of equilibrium $\widehat\varphi\equiv (\widehat{T},
\widehat{T^{*}},\widehat{V})$.
Then the stationary solution $\widehat\varphi$
is locally asymptotically stable (in $M_F$).
\end{theorem}

 \begin{remark}\label{sdd-vir3-eta-2}
 Similar to ODE case, described in
 \cite[Remark 13]{Rezounenko-DCDS-B-2017}, we have the following property.
For any $u\in C^1([-h,b);X)$
one has for $t\in [0,b)$
$${d\over dt}\eta(u_t) = [(D \eta)(u_t)](\dot u_t),$$
where $[(D \eta)(u_t)](\cdot )$ is the Fr\'{e}chet derivative of
$\eta$ at point $u_t$. Hence, (for a solution in
$\varepsilon$-neighborhood of the stationary solution
$\widehat{\varphi}$) the estimate $|{d\over dt}\eta(u_t)| \le || (D
\eta)(u_t)||_{L(C;R)} \cdot ||\dot u_t ||_{C} \le  \varepsilon\, ||
(D
\eta)(u_t)||_{L(C;R)}$ guarantees the property 
\begin{equation}\label{sdd-vir3-eta}
\left| {d\over dt}\eta(u_t) \right| \le \alpha_\varepsilon, \text{ with }
 \alpha_\varepsilon\to 0, \text{ as } \varepsilon\to 0.
\end{equation}
due to the boundedness of  $|| (D \eta)(\psi)||_{L(C;R)}$ as
$\varepsilon\to 0$ (here $||\psi - \widehat{\varphi}||_C <
\varepsilon$).

  \end{remark}


{\it Proof of Theorem \ref{vir3-theorem-stab1}}.
Let us consider (point-wise) the following  auxiliary  functional
$$U^\mathrm{sdd-x}(t,x) \equiv  \left( T(t,x) - \widehat T -
\int^{T(t,x)}_{\widehat T} {f(\widehat{T},\widehat{V})\over
f(\theta,\widehat{V})} \, d\theta  \right)e^{-\omega h}
 + \widehat{T^{*}}\cdot v\left({T^{*}(t,x)\over \widehat{T^{*}}}\right)
$$
\begin{equation}\label{sdd-vir3-Lyapuniv-functional-1}
+ {\widehat V  \over N}\cdot
v\left({V(t,x)\over \widehat V}\right)
+ \delta  \widehat{T^{*}} \int^t_{t-\eta(u_t)}
v\left({f(T(\theta,x),V(\theta,x))\over f(\widehat{T},\widehat{V})}
\right) \, d\theta.
\end{equation}

Now we can introduce the following Lyapunov functional  with
state-dependent delay along a solution of (\ref{sdd-vir3-04})

\begin{equation}\label{sdd-vir3-Lyapuniv-functional-2}
U^\mathrm{sdd}(t) \equiv   \int_\Omega U^\mathrm{sdd-x}(t,x) \, dx.
\end{equation}

The form of the functional is standard except the low limit
of the last integral in (\ref{sdd-vir3-Lyapuniv-functional-1})
which is state-dependent.
This state-dependence was first considered in
\cite{Rezounenko-DCDS-B-2017} (see also \cite{Rezounenko-EJQTDE-2016}).
For the constant delay case, see e.g. \cite{McCluskey-Yang-NA-2015}.

Now, for the simplicity of presentation, we consider the point-wise
time derivative of the functional $U^\mathrm{sdd-x}(t,x)$ defined in
(\ref{sdd-vir3-Lyapuniv-functional-1}). This time derivative is
considered along classical solutions of (\ref{sdd-vir3-04}).
It gives the possibility to consider ${\partial T(t,x)\over \partial t},
{\partial \widehat{T^{*}}(t,x)\over \partial t}, {\partial V(t,x)\over \partial t}, $
for any $t>0$. The computations below are in a sense close to the
ones in \cite{McCluskey-Yang-NA-2015}, but here we have two additional
diffusion terms and the state-dependence in both the system
(\ref{sdd-vir3-04}) and the Lyapunov functional.
First we consider the integral term
$$ {\partial \over \partial t}
\left[\int^t_{t-\eta(u_t)}
v\left({f(T(\theta,x),V(\theta,x))\over f(\widehat{T},\widehat{V})}
\right) \, d\theta
\right]
$$
$$= v\left({f(T(t,x),V(t,x))\over f(\widehat{T},\widehat{V})}
\right) - v\left({f(T(t-\eta(u_t),x),V(t-\eta(u_t),x))\over f(\widehat{T},\widehat{V})}
\right) \left( 1 - {d \over d t} \eta(u_t)\right)
$$
$$
= v\left({f(T(t,x),V(t,x))\over f(\widehat{T},\widehat{V})}
\right) - v\left({f(T(t-\eta(u_t),x),V(t-\eta(u_t),x))\over f(\widehat{T},\widehat{V})}
\right) + S^\mathrm{sdd}(t,x),
$$
where we denoted for short
\begin{equation}\label{vir3-s-sdd}
S^\mathrm{sdd}(t,x) \equiv
v\left({f(T(t-\eta(u_t),x),V(t-\eta(u_t),x))\over f(\widehat{T},\widehat{V})}
\right) \cdot {d \over d t} \eta(u_t).
\end{equation}
 \begin{remark}\label{sdd-vir3-S-sdd}
The term $S^\mathrm{sdd}$ appears due to the presence of the state-dependent delay. It makes the technical calculations more
challenging. The sign of $S^\mathrm{sdd}$ is undefined, so we propose below
(see also \cite{Rezounenko-DCDS-B-2017,Rezounenko-EJQTDE-2016}) a way to compensate/bound $S^\mathrm{sdd}$ by other positive defined terms in ${\partial U^\mathrm{sdd-x} \over \partial t}$ to have the time derivative of the Lyapunov functional (along a solution)  negative defined relative to the equilibrium.
   \end{remark}

Now we differentiate
$${\partial U^\mathrm{sdd-x}(t,x) \over \partial t} =
\left(1-  {f(\widehat{T},\widehat{V})\over
f(T(t,x),\widehat{V})}   \right) e^{-\omega h}\cdot {\partial T(t,x)\over \partial t} +
\left( 1- {\widehat{T^{*}}\over {T^{*}(t,x)}}\right)  \cdot {\partial \widehat{T^{*}}(t,x)\over \partial t}
$$
$$+ {1  \over N}\cdot
\left( 1 - {\widehat V\over V(t,x)}\right)
\cdot {\partial V(t,x)\over \partial t}
+ \delta  \widehat{T^{*}}  v\left({f(T(t,x),V(t,x))\over f(\widehat{T},\widehat{V})}
\right)
$$
$$
 - \delta  \widehat{T^{*}}  v\left({f(T(t-\eta(u_t),x),V(t-\eta(u_t),x))\over f(\widehat{T},\widehat{V})}
\right) + \delta  \widehat{T^{*}} S^\mathrm{sdd}(t,x).
$$
$$ = \left(1-  {f(\widehat{T},\widehat{V})\over
f(T(t,x),\widehat{V})}   \right) e^{-\omega h}\cdot
\left( \lambda - d T(t,x) - f(T(t,x),V(t,x)) + d^1 \Delta T(t,x) \right) $$
$${+}\left( 1{-}{\widehat{T^{*}}\over {T^{*}(t,x)}}\right)
\cdot \left(  e^{-\omega h} f(T(t{-}\eta(u_t),x),V (t{-}\eta(u_t),x)){-}\delta T^{*}(t,x){+}d^2
\Delta T^{*}(t,x)\right)
$$
$$+ {1  \over N}\cdot
\left( 1 - {\widehat V\over V(t,x)}\right)
\cdot \left( N\delta T^{*}(t,x) - cV (t,x) + d^3 \Delta V(t,x)\right) 
$$
$$
+ \delta  \widehat{T^{*}}  v\left({f(T(t,x),V(t,x))\over f(\widehat{T},\widehat{V})}
\right) - \delta  \widehat{T^{*}}
v\left({f(T(t-\eta(u_t),x),V(t-\eta(u_t),x))\over f(\widehat{T},\widehat{V})}
\right)
$$
$$+ \delta  \widehat{T^{*}} S^\mathrm{sdd}(t,x).
$$
 Calculations, using (\ref{sdd-vir3-stationary-1}), particularly,
 $\lambda = d \widehat{T} + f(\widehat{T},\widehat{V})$ give
$${\partial U^\mathrm{sdd-x}(t,x) \over \partial t}
= d \cdot \widehat{T}  \left( 1-  {T(t,x)\over \widehat{T}}\right)
 \left(1-  {f(\widehat{T},\widehat{V})\over
f(T(t,x),\widehat{V})}   \right) e^{-\omega h}$$
$$  +
\left(1-  {f(\widehat{T},\widehat{V})\over
f(T(t,x),\widehat{V})}   \right)  e^{-\omega h} \cdot d^1 \Delta T(t,x)
+ \left( 1- {\widehat{T^{*}}\over {T^{*}(t,x)}}\right)  \cdot d^2 \Delta T^{*}(t,x)
$$
$$
+ {1  \over N}\cdot
\left( 1 - {\widehat V\over V(t,x)}\right)
\cdot  d^3 \Delta V(t,x) + f(\widehat{T},\widehat{V})e^{-\omega h} \cdot C^1
+ \delta  \widehat{T^{*}}  v\left({f(T(t,x),V(t,x))\over f(\widehat{T},\widehat{V})}
\right)
$$
\begin{equation}\label{sdd-vir3-10}
 - \delta  \widehat{T^{*}}  v\left({f(T(t-\eta(u_t),x),V(t-\eta(u_t),x))\over f(\widehat{T},\widehat{V})}
\right) + \delta  \widehat{T^{*}} S^\mathrm{sdd}(t,x).
\end{equation}
where, for short, we collected some terms as  $ C^1$. It is written  as follows
$$ C^1 = C^1(t,x)\equiv \left( 1- {f(\widehat{T},\widehat{V})\over
f(T(t,x),\widehat{V})}\right) \left( 1- {f(T(t,x),V(t,x))\over
f(\widehat{T},\widehat{V}) }\right) $$
$$+
\left( 1- {\widehat{T^{*}}\over T^{*}(t,x)}\right) \left(  {f(T(t-\eta(u_t),x),V (t-\eta(u_t),x))\over f(\widehat{T},\widehat{V})
} - {  T^{*}(t,x)\over \widehat{T^{*}}}\right)
$$
$$+
\left( 1- {\widehat{V^{*}}\over V(t,x)}\right) \left(   {  T^{*}(t,x)\over \widehat{T^{*}}} - {V(t,x)\over \widehat{V }} \right).
$$
Calculations show that
$$C^1= 3 + {f(T(t,x),V(t,x))\over f(T(t,x),\widehat{V}) } +  {f(T(t-\eta(u_t),x),V (t-\eta(u_t),x))\over f(\widehat{T},\widehat{V})}
$$
$$ -  {f(\widehat{T},\widehat{V})\over f(T(t,x),\widehat{V}) }- {f(T(t,x),V(t,x))\over f(\widehat{T},\widehat{V})) } -  {f(T(t-\eta(u_t),x),V (t-\eta(u_t),x))\cdot \widehat{T^{*}} \over f(\widehat{T},\widehat{V})
 \cdot T^{*}(t,x)}
$$
$$ -  {  T^{*}(t,x)\cdot  \widehat{V }\over \widehat{T^{*}}\cdot V(t,x)} -   {V(t,x)\over \widehat{V }} .
$$
In the above expression we see two positive and five negative fraction terms, so we write $3=-2+5$ and add the following zero term ($0 = \ln 1 $):
$$ 0= \ln
\left[ \left(
{f(T(t,x),V(t,x))\over f(T(t,x),\widehat{V}) } \cdot  {f(T(t-\eta(u_t),x),V (t-\eta(u_t),x))\over f(\widehat{T},\widehat{V})}
\right)^{-1}
 \cdot {f(\widehat{T},\widehat{V})\over f(T(t,x),\widehat{V}) } \times \right.
$$
$$ \left. \times {f(T(t,x),V(t,x))\over f(\widehat{T},\widehat{V})) }\cdot  {f(T(t-\eta(u_t),x),V (t-\eta(u_t),x))\cdot \widehat{T^{*}} \over f(\widehat{T},\widehat{V})
 \cdot T^{*}(t,x)} \cdot {  T^{*}(t,x)\cdot  \widehat{V }\over \widehat{T^{*}}\cdot V(t,x)} \cdot {V(t,x)\over \widehat{V }} \right],
$$
which is split on the sum of seven logarithms
to write shortly, using the Volterra function $v$
$$C^1 = v \left( {f(T(t,x),V(t,x))\over f(T(t,x),\widehat{V}) }\right) + v \left(  {f(T(t-\eta(u_t),x),V (t-\eta(u_t),x))\over f(\widehat{T},\widehat{V})}\right)
$$
$$ - v \left( {f(\widehat{T},\widehat{V})\over f(T(t,x),\widehat{V}) }\right)
- v \left( {f(T(t,x),V(t,x))\over f(\widehat{T},\widehat{V})) }\right)
$$
\begin{equation}\label{sdd-vir3-c1}
- v \left(  {f(T(t-\eta(u_t),x),V (t-\eta(u_t),x))\cdot \widehat{T^{*}} \over f(\widehat{T},\widehat{V})
 \cdot T^{*}(t,x)}\right)
 - v \left(   {  T^{*}(t,x)\cdot  \widehat{V }\over \widehat{T^{*}}\cdot V(t,x)} \right)
- v  \left(   {V(t,x)\over \widehat{V }} \right).
\end{equation}
As before, we use $v(s)=s-1-\ln\, s $.


Now we discuss the diffusion terms (the ones with  coefficients $d^i$)
in (\ref{sdd-vir3-10}). More precisely, we are interested in
the sign of these terms after integration by $x$ in $\Omega$.
Denote them, for short, as
$$D^\mathrm{diff-3}(t,x)\equiv \left(1-  {f(\widehat{T},\widehat{V})\over
f(T(t,x),\widehat{V})}   \right) e^{-\omega h}\cdot d^1 \Delta T(t,x)
 + \left( 1- {\widehat{T^{*}}\over {T^{*}(t,x)}}\right)  \cdot d^2 \Delta T^{*}(t,x)
$$
\begin{equation}\label{sdd-vir3-Ddiff-3}
+ {1  \over N}\cdot
\left( 1 - {\widehat V\over V(t,x)}\right)
\cdot  d^3 \Delta V(t,x), \qquad D^\mathrm{diff-3}(t)\equiv \int_{\Omega} D^\mathrm{diff-3}(t,x)\, dx.
\end{equation}

In case of differentiable $f$  $(Hf_4)$ (see (\ref{sdd-vir3-f-3})) need the following simple
\begin{proposition}\label{proposition-A}

Let $p: \mathbb{R} \to \mathbb{R}$ be differentiable.  
Then $\int_\Omega p(u(x)) \Delta u(x)\, dx = - \int_\Omega p^\prime (u) ||\nabla u||^2 \, dx$ 
for any
$u\in C^2(\overline{\Omega})$ satisfying ${\partial u(x)\over \partial n} |_{\partial \Omega}=0$.
\end{proposition}

{\it Proof of Proposition \ref{proposition-A}.} We use the classical
Gauss-Ostrogradsky theorem. Consider the vector field $E\equiv p(u) \nabla u$.
Hence $\mathrm{div}\, E = p^\prime (u) ||\nabla u||^2 + p(u) \Delta u$.
One has $\int_\Omega \mathrm{div}\, E \, dx = \int_\Omega p^\prime (u)
||\nabla u||^2 \, dx + \int_\Omega p(u) \Delta u \, dx  =
\int_{\partial\Omega} p(u) ( \nabla u, n)\, dS =0.$ The last
equality due to the Neumann boundary conditions. Finally,
$\int_\Omega p(u) \Delta u \, dx  = - \int_\Omega p^\prime (u) ||\nabla u||^2 \, dx.$
It completes the proof of Proposition \ref{proposition-A}.
\smallskip


Now we apply Proposition \ref{proposition-A} to show that $D^\mathrm{diff-3}(t)\le 0.$ Let us start with the first term in $D^\mathrm{diff-3}(t,x)$, see (\ref{sdd-vir3-Ddiff-3}), and show that
$  \int_\Omega \left( 1-{f(\widehat{T},\widehat{V})\over f(T(t,x),\widehat{V})}\right) \, \Delta T(t,x)\, dx \le 0.$ For this we set $p(T)=\left( 1-{f(\widehat{T},\widehat{V})\over f(T(t,x),\widehat{V})}\right)$ and check that
$p^\prime (T) =
f^\prime_1(T,\widehat{V}) f(\widehat{T},\widehat{V})\times [ f(T,\widehat{V})]^{-2} \ge 0$ due to $f^\prime_1(T,\cdot) \ge 0$ by the assumption on $f$. Similar considerations with the second and third terms in (\ref{sdd-vir3-Ddiff-3}) show that
$$ D^\mathrm{diff-3}(t)\equiv \int_{\Omega} D^\mathrm{diff-3}(t,x)\, dx
$$
$$=
 -  d^1 \cdot e^{-\omega h}f(\widehat{T},\widehat{V})\int_{\Omega} {f^\prime_1(T(t,x),\widehat{V}) \over f(T(t,x),\widehat{V})^2} \, ||\nabla T(t,x)||^2 \, dx
$$
\begin{equation}\label{sdd-vir3-D1-4}
- d^2 \cdot\widehat{T^{*}} \int_{\Omega}
{||\nabla T^{*}(t,x)||^2 \over [T^{*}(t,x)]^2} \, dx
- d^3 {\widehat V  \over N}\cdot \int_{\Omega}
 { ||\nabla V(t,x)||^2\over [V(t,x)]^2}
\cdot
 \, dx\le 0.
\end{equation}

\begin{remark}\label{remark-B}
In case of nondifferentiable $f$ we prove $D^\mathrm{diff-3}(t)\le 0$,
using alternative (geometrical) conditions on $f$ given
in $(Hf_4)$ (see (\ref{sdd-vir3-f-3})).

Then $D^\mathrm{diff-3}(t)\le 0$ along any classical solution.

\end{remark}


Now we combine the arguments above to study the Lyapunov functional
$U^\mathrm{sdd}(t)$, see (\ref{sdd-vir3-Lyapuniv-functional-2}). We have
the following equality (c.f. (\ref{sdd-vir3-10}))
$$ {d\over dt} U^\mathrm{sdd}(t) {=} \int_{\Omega}
{\partial U^\mathrm\mathrm(t,x) \over \partial t} \, dx
= d  \widehat{T} \cdot  e^{-\omega h}\int_{\Omega}
 \left( 1{-}{T(t,x)\over \widehat{T}}\right)
 \left(1{-}{f(\widehat{T},\widehat{V})\over
f(T(t,x),\widehat{V})}   \right)    \, dx
$$
 $$ + D^\mathrm{diff-3}(t) + f(\widehat{T},\widehat{V})e^{-\omega h}
  \cdot\int_{\Omega}  C^1  \, dx
$$
$$+ \delta  \widehat{T^{*}}\int_{\Omega}
\left[   v\left({f(T(t,x),V(t,x))\over f(\widehat{T},\widehat{V})}
\right) -   v\left({f(T(t-\eta(u_t),x),V(t-\eta(u_t),x))\over f(\widehat{T},\widehat{V})}
\right) \right.
$$
$$\left. + \, S^\mathrm{sdd}(t,x) \rule{0pt}{5mm}  \right]  \, dx.
 $$
 Here $D^\mathrm{diff-3}(t)$ is defined in (\ref{sdd-vir3-Ddiff-3}) and transformed in (\ref{sdd-vir3-D1-4}), $C^1$ is presented as in (\ref{sdd-vir3-c1}) and $S^\mathrm{sdd}$ is defined in (\ref{vir3-s-sdd}).
 We remind (see (\ref{sdd-vir3-stationary-1})) that $\delta  \widehat{T^{*}} = e^{-\omega h} f(\widehat{T},\widehat{V}) $  which leads  to cancellation of the first and second terms in the last integral with the corresponding terms in $C^1$ (see (\ref{sdd-vir3-c1})). We continue calculations
 $$ {d\over dt} U^\mathrm{sdd}(t) =
\int_{\Omega} {\partial U^\mathrm{sdd-x}(t,x) \over \partial t} \, dx =
 d  \widehat{T} \cdot  e^{-\omega h}\int_{\Omega}  \left( 1-  {T(t,x)\over \widehat{T}}\right) \left(1-  {f(\widehat{T},\widehat{V})\over
f(T(t,x),\widehat{V})}   \right)    \, dx
$$
 $$ + f(\widehat{T},\widehat{V})e^{-\omega h} \cdot\int_{\Omega}
\left\{ - v \left( {f(\widehat{T},\widehat{V})\over f(T(t,x),\widehat{V}) }\right)
{-} v \left(  {f(T(t-\eta(u_t),x),V (t-\eta(u_t),x))\cdot \widehat{T^{*}} \over f(\widehat{T},\widehat{V})
 \cdot T^{*}(t,x)}\right)
\right. dx
$$
$$ \left.- v \left(   {  T^{*}(t,x)\cdot  \widehat{V }\over \widehat{T^{*}}\cdot V(t,x)} \right)
- \left[ v  \left(   {V(t,x)\over \widehat{V }} \right) - v \left( {f(T(t,x),V(t,x))\over f(T(t,x),\widehat{V}) }\right)\right]
\right\} \, dx
$$
\begin{equation}\label{sdd-vir3-11}
+ D^\mathrm{diff-3}(t) + \delta  \widehat{T^{*}}\int_{\Omega}  S^\mathrm{sdd}(t,x) \, dx.
\end{equation}
We will show that  all the terms in (\ref{sdd-vir3-11}) are non-positive except for the last one
which, in general, may change sign. The first term in (\ref{sdd-vir3-11}) is non-positive due to monotonicity of $f$ with respect to the first coordinate. The property $D^\mathrm{diff-3}(t)\le 0 $ is given in (\ref{sdd-vir3-D1-4}).
To show that
$$
\int_{\Omega} \left[ v  \left(   {V(t,x)\over \widehat{V }} \right) - v \left( {f(T(t,x),V(t,x))\over f(T(t,x),\widehat{V}) }\right)\right]
 \, dx \ge 0
$$
 we use the property  $(Hf_3)$ of function $f$ (see (\ref{sdd-vir3-12})).

Now we plan to prove that  ${d\over dt} U^\mathrm{sdd}(t) \le 0$ in a
small neighbourhood of the stationary solution with the equality
only in case of $(T,T^{*},V)=(\widehat{T },\widehat{T^{*}},\widehat{V})$.
In the particular case of constant delay, one has $S^\mathrm{sdd}(t,x)=0$
which may lead to the global stability of $(\widehat{T },\widehat{T^{*}},\widehat{V})$.

We rewrite, for short, (\ref{sdd-vir3-11}) as
\begin{equation}\label{sdd-vir3-13}
{d\over dt} U^\mathrm{sdd}(t) =  \delta  \widehat{T^{*}}\int_{\Omega}
\left( \rule{0pt}{4mm}
 - D^\mathrm{sdd}(t,x) + S^\mathrm{sdd}(t,x) \right) \, dx,
\end{equation}
where $D^\mathrm{sdd}(t,x)$ contains all the terms except the last one in (\ref{sdd-vir3-11}).
As proved above $  \int_{\Omega} D^\mathrm{sdd}(t,x)  \, dx \ge 0$.
Let us start with an analysis of the zero-sets
 $D^\mathrm{sdd}(t,x)=0, S^\mathrm{sdd}(t,x)=0$ and ${d\over dt} U^\mathrm{sdd}(t) = 0$.

We start with $D^\mathrm{sdd}(t,x)=0$. One sees from (\ref{sdd-vir3-11}) that $T=\widehat{T}$.
Since $v(s)=0$ iff $s=1$, we see from (\ref{sdd-vir3-12}) that $V=\widehat{V}$.
 Hence $T^{*}=\widehat{T^{*}}$. One also sees
 $f(T(t-\eta(u_t),x),V (t-\eta(u_t),x))=f(\widehat{T},\widehat{V})$.
 Moreover, $D^\mathrm{diff-3}(t) =0$, means (see Proposition \ref{proposition-A} and (\ref{sdd-vir3-D1-4}))
  that $T, T^{*}$ and $V$ are independent of $x$.
The zero set $S^\mathrm{sdd}(t,x)=0$ is described (see (\ref{vir3-s-sdd})
  by  $f(T(t-\eta(u_t),x),V (t-\eta(u_t),x))=f(\widehat{T},\widehat{V})$
  or  ${d \over d t} \eta(u_t)=0$ along a solution.
It is important for us that the zero-set $D^\mathrm{sdd}(t,x)=0$ is a singleton
$(T,T^{*},V)=(\widehat{T},\widehat{T^{*}},\widehat{V})$ and is a subset
of $S^\mathrm{sdd}(t,x)=0$. The rest of the proof that in a small neighbourhood
of $(\widehat{T},\widehat{T^{*}},\widehat{V})$ one has
$|S^\mathrm{sdd}(t,x)| < D^\mathrm{sdd}(t,x)$ follows the streamline of the proof
\cite[Theorem 12]{Rezounenko-DCDS-B-2017} (see also \cite[Theorem 3.3]{Rezounenko-EJQTDE-2016}).
It relies on property (\ref{sdd-vir3-eta}), auxiliary quadratic functionals due to property (\ref{sdd-vir3-v-functional-2}) of Volterra function $v$ and the change of variables to the polar ones (see \cite[(33)-(35)]{Rezounenko-DCDS-B-2017}). We do not repeat the calculations here. The property ${d\over dt} U^\mathrm{sdd}(t) \le 0$ in a small neighbourhood of the stationary solution with the equality only in case of $(T,T^{*},V)=(\widehat{T },\widehat{T^{*}},\widehat{V})$ completes the proof of Theorem \ref{vir3-theorem-stab1}.

\medskip

It is interesting to notice that
 $\varphi \in M_F$ (see (\ref{sdd-vir3-M_F})) is not a necessary condition for our approach. Now we consider a wider set $\Omega_{Lip}$ (see (\ref{sdd-vir3-omega})).
 Let us discuss a particular simple form of the delay (c.f. examples in \cite{Rezounenko-DCDS-B-2017})
\begin{equation}\label{sdd-vir3-example1}
\eta(\varphi) = \int^{0}_{-h} \, \xi(\varphi(\theta))\,
d\theta, \qquad \varphi \in C
\end{equation}
 with a locally Lipschitz $\xi$.
To check the property (\ref{sdd-vir3-eta})  we
calculate
$$ {d\over dt} \eta(u_t) =  {d\over dt}
\int^{0}_{-h} \, \xi(u(t+\theta))\, d\theta =
 {d\over dt} \int^{t}_{t-h} \, \xi(u(s))\, ds = \xi(u(t)) -
 \xi(u(t-h)).
$$
Hence, in the $\varepsilon$-neighborhood of the stationary solution $\hat
u$, one has
$$ \left| {d\over dt} \eta(u_t) \right| \le  \left| \xi(u(t)) -
 \xi(u(t-h))\right| \le 2 \varepsilon\, L_{\xi,\varepsilon}\equiv \alpha_\varepsilon \to 0 \quad \hbox{ as } \varepsilon \to 0.
$$
Here $L_{\xi,\varepsilon}$ is the Lipschitz constant of $\xi$.
More general delay terms could be used
\begin{equation}\label{sdd-vir3-example2}
\eta(\varphi) = \rho \left(\int^{0}_{-h} \, \xi(\varphi(\theta)) \kappa(\theta)\,
d\theta \right), \qquad \varphi \in C, \quad \kappa\in
C([-h,0];R) %
\end{equation}
with a differentiable $\rho: \mathbb{R} \to [0,h]$. The example
(\ref{sdd-vir3-example1}) is a particular case of (\ref{sdd-vir3-example2}) with $\rho(s) \equiv
s$ and $\kappa(s)\equiv 1.$

The discussion above shows that property (\ref{sdd-vir3-eta})  of the state-dependent delay (\ref{sdd-vir3-example2}) allows to use the proof of Theorem \ref{vir3-theorem-stab1} to get  the following result in $\Omega_{Lip}$
\begin{theorem}\label{vir3-theorem-stab2}
Let non-linear function $f$ satisfy $(Hf_1+), (Hf_2), (Hf_3), (Hf_4)$
(see (\ref{sdd-vir3-f-1+}), (\ref{sdd-vir3-f-3}),
(\ref{sdd-vir3-12})) and  state-dependent delay $\eta : C \to [0,h]$ be of the form (\ref{sdd-vir3-example2}).
Then the stationary solution $\widehat\varphi$
is locally asymptotically stable.
\end{theorem}

% ...

{\bf Acknowledgments.} The author is thankful to anonymous referees for useful comments and suggestions. This work was supported in part by GA CR under project 16-06678S.



\medskip

Submitted February 27, 2017. 


\medskip


\begin{thebibliography}{99}




\bibitem{Beddington-JAE-1975} 
 J. R. Beddington,
  {Mutual interference between parasites or predators and its effect on searching efficiency},
 \emph{Journal of Animal Ecology,} \textbf{44} (1975), 331--340.



\bibitem{Carloni et al-CMM-2012} 
 G.Carloni, A.Crema, M.B. Valli, A.Ponzetto,  M.Clementi,
  {HCV Infection by Cell-to-Cell Transmission: Choice or Necessity?}
 \emph{Current Molecular Medicine}, \textbf{12} (2012), 83--95.





\bibitem{Chueshov-Rezounenko_CPAA-2015} 
 I.D. Chueshov, A.V. Rezounenko,
  {Finite-dimensional global attractors for parabolic nonlinear equations with state-dependent delay},
 \emph{Communications on Pure and Applied Analysis,} \textbf{14/5} (2015), 1685-1704.

\bibitem{DeAngelis-Goldstein-ONeill-E-1975} 
 D. L. DeAngelis, R. A. Goldstein and R. V. O'Neill,
  {A model for tropic interaction}, \emph{Ecology}, \textbf{56} (1975), 881--892.








\bibitem{Walther_book} 
 O. Diekmann, S. van Gils, S. Verduyn Lunel and H.-O. Walther,
 \emph{Delay Equations: Functional, Complex, and Nonlinear Analysis},
 Springer-Verlag, New York, 1995.


\bibitem{Driver-AP-1963} 
 R. D. Driver,
  {A two-body problem of classical electrodynamics: The one-dimensional case},
 \emph{Ann. Physics}, \textbf{21} (1963), 122--142.

\bibitem{Gourley-Kuang-Nagy-JBD-2008}
  S.A. Gourley, Y.Kuang, J.D. Nagy,
  {Dynamics of a delay differential equation model
of hepatitis B virus infection},
 \emph{Journal of Biological Dynamics} 2, (2008)  140--153.


\bibitem{Hale} 
 J. K. Hale,
 \emph{Theory of Functional Differential Equations},
 Springer, Berlin- Heidelberg- New York, 1977.





\bibitem{Hartung-Krisztin-Walther-Wu-2006}  
 F.~Hartung, T.~Krisztin, H.-O.~Walther and J.~Wu,
  {Functional differential equations with state-dependent delays: Theory and applications},
 In: \emph{Canada, A., Drabek., P. and A. Fonda (Eds.) Handbook of Differential Equations, Ordinary Differential Equations,} Elsevier Science B.V., North Holland, \textbf{3} (2006), 435--545.


\bibitem{Hattaf-Yousfi-CMA-2015} 
 K. Hattaf,  N. Yousfi,
  {A generalized HBV model with diffusion and two delays},
 \emph{Computers and Mathematics with Applications}, \textbf{69} (2015), 31–-40.




 \bibitem{Huang-Ma-Takeuchi-AML-2011} 
 G. Huang, W. Ma, Y. Takeuchi,
  {Global analysis for delay virus
dynamics model with Beddington-DeAngelis functional response},
 \emph{Applied Mathematics Letters}, 24 (2011) 1199--1203.


\bibitem{Korobeinikov-BMB-2007} 
 A. Korobeinikov,
  {Global properties of infectious disease models with nonlinear incidence},
 \emph{Bull. Math. Biol.}, \textbf{69} (2007), 1871--1886.


\bibitem{Kuang-1993_book} 
 Y. Kuang,
 \emph{Delay Differential Equations with Applications in Population Dynamics},
 Mathematics in Science and Engineering, 191. Academic Press, Inc., Boston, MA, 1993.


\bibitem{Lyapunov-1892}
 A. M. Lyapunov,
 \emph{The General Problem of the Stability of Motion},
 Kharkov Mathematical Society, Kharkov, 1892, 251p.



\bibitem{Mallet-Paret} 
 J. Mallet-Paret,  R. D. Nussbaum, P. Paraskevopoulos,
 {Periodic solutions for functional-differential equations  with multiple state-dependent time lags},
 \emph{Topol. Methods Nonlinear Anal.}, \textbf{3:1} (1994), 101--162.






\bibitem{Martin-Smith-TAMS-1990} 
 R.H. Martin, Jr., H.L. Smith,
  { Abstract functional-differential equations and reaction-diffusion systems},
 \emph{Trans. Amer. Math. Soc.}, \textbf{321} (1990), 1--44.


%


\bibitem{McCluskey-Yang-NA-2015} 
 C. McCluskey,    Yu.Yang,
  {Global stability of a diffusive virus dynamics model with general incidence function and time delay},
 \emph{Nonlinear Anal. Real World Appl}, \textbf{25} (2015), 64-78.


\bibitem{Murray-Kelleher-Cooper_J. Virol-2011} 
  J.M. Murray, A.D. Kelleher, D.A. Cooper,
  {Timing of the Components of the
HIV Life Cycle in Productively Infected CD4+ T Cells in a Population
of HIV-Infected Individuals},
 \emph{J. Virol.} (2011), vol. 85 no.
20, 10798-10805.




\bibitem{Nowak-Bangham-S-1996} 
 M. Nowak and C. Bangham,
  {Population dynamics of immune response to persistent viruses},
 \emph{Science}, \textbf{272} (1996), 74--79.


\bibitem{Pawlotsky-Semin Liver Dis-2014} 
 JM. Pawlotsky,
  {New hepatitis C virus (HCV) drugs and the hope for a cure: concepts in anti-HCV drug development},
 \emph{Semin Liver
Dis.}, \textbf{34(01)} (2014), 22--29.



\bibitem{Pazy-1983-book} 
 A. Pazy,
  {\emph{Semigroups of linear operators and applications to partial differential equations. Applied Mathematical Sciences}},
 44. Springer-Verlag, New York, 1983. {\rm viii}+279 pp.




\bibitem{Perelson-Neumann-Markowitz-Leonard-Ho-S-1996} 
 A. Perelson, A. Neumann, M. Markowitz, J. Leonard and D. Ho,
  {HIV-1 dynamics in vivo: Virion clearance rate, infected cell life-span, and viral generation time},
 \emph{Science}, \textbf{271} (1996), 1582--1586.



\bibitem{Rezounenko_JMAA-2007} 
 A. V. Rezounenko,
  {Partial differential
equations with discrete and distributed state-dependent delays},
 \emph{Journal of Mathematical Analysis and Applications}, \textbf{326} (2007), 1031--1045.


\bibitem{Rezounenko_NA-2009} 
 A. V. Rezounenko,
  {Differential equations with discrete state-dependent delay: Uniqueness and well-posedness in the space of continuous functions},
 \emph{Nonlinear Analysis: Theory, Methods and Applications}, \textbf{70} (2009), 3978--3986.


\bibitem{Rezounenko_NA-2010} 
 A. V. Rezounenko,
  {Non-linear partial differential equations with discrete state-dependent delays in a metric space},
 \emph{Nonlinear Analysis: Theory, Methods and Applications}, \textbf{73} (2010), 1707--1714.




\bibitem{Rezounenko_JMAA-2012} 
 A. V. Rezounenko,
  {A condition on delay for differential equations with discrete state-dependent delay},
 \emph{Journal of Mathematical Analysis and Applications}, \textbf{385} (2012), 506--516.

\bibitem{Rezounenko-JADEA-2012} 
 A. V. Rezounenko,
 Local properties of solutions to non-autonomous parabolic PDEs with state-dependent delays,
 \emph{Journal of Abstract Differential Equations and Applications}, \textbf{2} (2012), 56--71.



\bibitem{Rezounenko-Zagalak-DCDS-2013} 
 A.V. Rezounenko,  P. Zagalak,
  {Non-local PDEs
with discrete state-dependent delays: well-posedness in a metric space},
 \emph{Discrete and Continuous Dynamical Systems - Series A}, \textbf{33:2} (2013), 819--835.



\bibitem{Rezounenko-DCDS-B-2017} 
 A. V. Rezounenko,
  {Stability of a viral infection model with state-dependent delay, CTL and antibody immune responses},
 \emph{Discrete and Continuous Dynamical Systems
- Series B}, \textbf{22} (2017), 1547--1563; Preprint arXiv:1603.06281v1 [math.DS], 20 March 2016,
\href{http://arxiv.org/abs/1603.06281v1}{arxiv.org/abs/1603.06281v1}.


\bibitem{Rezounenko-EJQTDE-2016} 
 A. V. Rezounenko,
  {Continuous solutions to a viral infection model with general incidence rate, discrete state-dependent delay, CTL and antibody immune responses},
 \emph{Electron. J. Qual. Theory Differ. Equ.}, \textbf{79} (2016), 1--15.










\bibitem{Shudo-Ribeiro-Talal-Perelson_Antiviral
Therapy-2008} 
 E. Shudo,  R.M. Ribeiro,  A.H. Talal,  A.S.
Perelson,
  {A hepatitis C viral kinetic model that allows for time-varying drug effectiveness},
 \emph{Antiviral Therapy}, \textbf{13} (2008), 919--926.






\bibitem{Smith-1995-book} 
 H. L. Smith,
 \emph{Monotone Dynamical Systems. An Introduction to the Theory of Competitive and Cooperative Systems},
 Mathematical Surveys and Monographs, 41. American Mathematical Society, Providence, RI, 1995.


\bibitem{Smith-2011-book}
  H. Smith,
 An
Introduction to Delay Differential Equations with Sciences
Applications to the Life, Texts in Applied Mathematics, vol. 57,
Springer, New York, Dordrecht, Heidelberg, London, 2011. 





\bibitem{travis_webb} 
 C.C. Travis and G.F. Webb,
  {Existence and
stability for partial functional differential equations},
 \emph{Transactions of AMS}, \textbf{200} (1974), 395--418.





\bibitem{Walther_JDE-2003} 
 H.-O. Walther,
  {The solution manifold and $C\sp 1$-smoothness for differential equations with state-dependent delay},
 \emph{Journal of Differential Equations}, \textbf{195} (2003), 46--65.


\bibitem{Wang-Liu-MMAS-2013} 
  X. Wang, S. Liu,
  A class of
delayed viral models with saturation infection rate and immune
response,
  \emph{Math. Methods Appl. Sci.} 36,  \textbf{2} (2013),
125--142. 




\bibitem{Wang-Huang-Zou-AA-2014} 
 F.-B. Wang, Y. Huang, X. Zou,
  {Global dynamics of a PDE
in-host viral model},
 \emph{Applicable Analysis: An International Journal}, \textbf{93:11} (2014), 2312--2329.



\bibitem{Wang-Wang-MB-2007_HBV_spatial dependence} 
 K. Wang, W. Wang,
  {Propagation of HBV with spatial dependence},
 \emph{ Math. Biosci.}, \textbf{201} (2007), 78--95.


\bibitem{Wang-Yang-Kuniya-JMAA-2016} 
 J.Wang, J.Yang, T.Kuniya,
   {Dynamics of a PDE viral infection model incorporating cell-to-cell transmission},
 Journal of Mathematical Analysis and Applications,
444, (2016), 1542-1564.





\bibitem{WHO-Global hepatitis report-2017}
 World Health Organization,
\emph{Global hepatitis report-2017, April 2017,} ISBN: 978-92-4-156545-5 \href{http://apps.who.int/iris/bitstream/10665/255016/1/9789241565455-eng.pdf?ua=1}{http://apps.who.int/iris/bitstream/10665/255016/1/9789241565455-eng.pdf?ua=1}


\bibitem{Wu_book} 
 J. Wu,
 \emph{Theory and Applications of Partial Functional
Differential Equations},
 Springer-Verlag, New York, 1996.



\bibitem{Xu_JQTDE-2012} 
 S. Xu,
  {Global stability of the
virus dynamics model with Crowley-Martin functional response},
 \emph{J. Qual. Theory Differ. Equ.}, \textbf{2012}(9),  (2012), 1--10.










\bibitem{Zhao-Xu_EJDE-2014} 
 Y. Zhao, Z. Xu,
 Global dynamics for a delayed hepatitis C virus infection model,
 \emph{Electronic Journal of Differential Equations}, \textbf{2014}/132  (2014), 1--18.





\end{thebibliography}
\end{document}